\documentclass[twoside,reqno,11pt]{fcaa-var} %
\usepackage{graphicx, epsfig, latexsym}
\usepackage{amsthm, amsmath, amsfonts, amssymb}

\newtheoremstyle{theorem}
  {15pt}          
  {15pt}  
  {\sl}  
  {\parindent}
  {\sc}  
  {. }   
  { }    
  {}     
\theoremstyle{theorem}

\newtheoremstyle{defi}
  {15pt}          
  {15pt}  
  {\rm}  
  {\parindent}     
  {\sc}  
  {. }    
  { }    
  {}     
\theoremstyle{defi}

 \def\proofname{\indent\rm P\ r\ o\ o\ f.}
 
 \def\qed{\hfill$\Box$} 

 \usepackage{hyperref} 

 \def\theequation{\arabic{section}.\arabic{equation}}

 \def\themycopyrightfootnote{ 
   }

 \title[Multivariate Grunwald-Letnikov]{Infinitesimal translations and a multivariate Gr\"unwald-Letnikov calclulus}
 \author[\normalsize A. Pallavi Sudhir]{\normalsize Abhimanyu Pallavi Sudhir $^1$}

\begin{document}



 \bigskip \medskip

 \begin{abstract}

The goal of this paper is to construct a multivariate generalisation of the Gr\"unwald-Letnikov derivative, a classical fractional derivative operator. To do so, we first produce a formalism of fractional derivatives in terms of infinitesimal translations that justifies the ``binomial theorem'' argument for the Gr\"unwald-Letnikov derivative, allowing us to then extend the argument to construct the multivariate derivative via the more general multi-binomial theorem. We conclude by studying the principal value of the fractional derivative of a multivariate power function, obtaining a characteristic equation in agreement with recent research in the area.
\medskip


 \smallskip


 \end{abstract}

 \maketitle

 \vspace*{-16pt}


 \section{Introduction}\label{sec:1}

\setcounter{section}{1}
\setcounter{equation}{0}\setcounter{theorem}{0}

The Gr\"unwald-Letnikov derivative is a generalisation of the higher-order derivative to fractional orders -- the classic motivation for the derivative comes from considering the limit form of the $n$\textsuperscript{th}-order derivative, observing the similarity to the binomial theorem and writing down an analogous generalisation to the binomial series \cite{diaz1974}.

\begin{equation}
  D^Rf(x) = \lim\limits_{h\to0}h^{-R}\sum\limits_{k=0}^\infty \binom{R}{k}(-1)^kf(x-kh)
\label{gl-def}
\end{equation}

This paper makes a similar observation with the multi-binomial theorem to define the multivariate Gr\"unwald-Letnikov derivative -- however, to rigorously justify doing so, we will first introduce a formalism of ``infinitesimal function translations'' that will make precise the relation between the $n$\textsuperscript{th}-order derivative and the Binomial theorem.

As part of our exploration of this multivariable Gr\"unwald-Letnikov derivative, we generalise the results of \cite{pasu2018} to a multivariable setting. In \cite{pasu2018}, a relation between terms $h$ and $N$ of the limits $h \to 0$ and the implicit limit of summation $N \to \infty$ was enforced on the Gr\"unwald-Letnikov derivative to ``force'' it to equal its ``principal value'', which is definitionally the Riemann-Liouville derivative. The condition is of the form $h = qx/N$, where the parameter $q$ depends on the function and the order of the derivative -- for the simple case of the power function, $D^R x^p$, it is known that $q$ is the solution to the following equation, known as the ``characteristic equation'' of the derivative: \cite{pasu2018}

\begin{equation}
q^{-R} {}_2F_1(-p, -R; 1 - R; q) = \frac{\pi R}{\sin \pi R}\binom{p}{R}
\label{charpol}
\end{equation}

Using the multivariate Riemann-Liouville derivative defined in \cite{malk2015}, we derive a generalisation of Eq.~\eqref{charpol} to a multivariable power function of the form ${x_1}^{p_1}\ldots{x_n}^{p_n}$ -- as we will see, this will allow us to write down a value of $q$ for any analytic function.


\section{The infinitesimal translation formalism}\label{sec:2}

\setcounter{section}{2}
\setcounter{equation}{0}\setcounter{theorem}{0}

Consider introducing a translation operator $\phi^h$ that translated a function $f : \mathbb{R} \to \mathbb{R}$ left-ward by a real number $h$, i.e.

\begin{equation}
  \phi ^ h f(x) \triangleq f(x + h)
\label{gen-def}
\end{equation}

It can be verified that such operators generate (with addition and scalar multiplication) a commutative ring under extensional addition and operator multiplication with $1 := \phi^0$. Further, it is clear that the superscript operation $\phi^h$ is equal to the integer power for integer $h$. One may then write the derivative as:

\begin{equation}
  \frac{d}{dx} = \lim\limits_{h \to 0} \frac{\phi ^ h - 1}{h}
\label{der-def}
\end{equation}

Working with the $n$\textsuperscript{th}-order derivative $d^n/dx^n$ is then straightforward. The binomial theorem is true on any commutative ring, therefore we may write, where the limit exists:

\begin{equation}
  \frac{d^n}{dx^n} = \lim\limits_{h \to 0} \frac{1}{h ^ n} \sum_{k = 0}^\infty (-1)^k \binom{n}{k} \phi^{(n - k)h}
\label{phi-binom}
\end{equation}

The Gr\"unwald-Letnikov derivative then simply becomes the formalisation of a real power ($n \in \mathbb{R}$) on the ring of translation operators via the binomial expansion as can be seen in Eq.~\eqref{phi-binom} by writing $\phi^{(n - k)h}f = f(x + (n - k)h)$.

\section{The multivariate fractional derivative}\label{sec:3}

\setcounter{section}{3}
\setcounter{equation}{0}\setcounter{theorem}{0}

With this formalism in mind, it is then much easier to motivate the multivariate Gr\"unwald-Letnikov derivative. For functions $f : \mathbb{R}^n \to \mathbb{R}^m$, we may define $n$ generators $\phi_i$:

\begin{equation}
  {\phi_i}^h f(x_1,\ldots x_i,\ldots x_n) \triangleq f(x_1,\ldots x_i + h,\ldots x_n)
\label{gens-def}
\end{equation}

And the partial derivative in the $x_i$ direction is:

\begin{equation}
  \frac{\partial}{\partial x_i} = \lim\limits_{h \to 0} \frac{{\phi_i} ^ h - 1}{h}
\label{ders-def}
\end{equation}

We are then interested in the general-order mixed partial derivative, which we write in terms of the infinitesimal translation operators as follows (the $h_i$'s can be set to each other should $f$ be well-behaved, implying symmetry of the mixed derivative):

\begin{equation}
  \frac{{{\partial ^{\sum r_i}}}}{{\partial {x_1}^{{r_1}} \ldots \partial {x_n}^{{r_n}}}} = \lim\limits_{h_i \to 0}{\left( {\frac{{{\phi _1}^{h_1} - 1}}{h_1}} \right)^{{r_1}}} \ldots {\left( {\frac{{{\phi _n}^{h_n} - 1}}{h_n}} \right)^{{r_n}}}
\label{hders-def}
\end{equation}

Analogous to the significance of the binomial theorem in the univariate case, we can see that the limit form of the general-order partial derivative will match the form of the \textit{multi-binomial theorem}, with which we expand the above expression (the limit is suppressed, and $i$ runs from 1 to $n$):

\begin{equation}
  \frac {\partial ^{\sum r_i}}
        {\partial {x_1}^{r_1} \ldots \partial {x_n}^{r_n}} = 
  \prod\limits_i {h_i}^{-r_i} \cdot 
  \sum  \limits_{k_i = 0}^{r_i}
        {\prod\limits_i (-1)^{k_i}\binom{r_i}{k_i}{\phi_i}^{-k_i h_i}} 
\label{phis-multibinom}
\end{equation}

Where we replaced ${\phi_i}^{(r_i-j_i)h}$ with ${\phi_i}^{-j_i h}$ as $r_i h$ is infinitesimal and can be factored out of the summation. Our definition of the multivariable Gr\"unwald-Letnikov derivative for non-integer $r_i$ is then simply obtained by replacing the upper limits of the summation with $\infty$. Or in shorthand where $\circ$ represents pointwise multiplication and $k_i = \mathbf{k} \cdot \mathbf{e}_i$, $h_i = \mathbf{h} \cdot \mathbf{e}_i$,

\begin{equation}
  \frac {\partial ^{\sum r_i} f(\mathbf{x})}
        {\partial {x_1}^{r_1} \ldots \partial {x_n}^{r_n}} = 
  \prod\limits_i {h_i}^{-r_i} \cdot 
  \sum  \limits_{\mathbf{k} \in \mathbb{N}^n} 
        {\left[ \prod\limits_i (-1)^{k_i}\binom{r_i}{k_i} \right]
  f(\mathbf{x} - \mathbf{k} \circ \mathbf{h})} 
\label{mvglder-def}
\end{equation}
\vspace*{2pt} 

\section{Considerations on the $q$-principal value}\label{sec:4}

\setcounter{section}{4}
\setcounter{equation}{0}\setcounter{theorem}{0}

Consider the derivative: $\frac {\partial ^{\sum r_i}} {\partial {x_1}^{r_1} \ldots \partial {x_n}^{r_n}} \left({x_1}^{p_1}\ldots{x_n}^{p_n}\right)$. The principal value of this derivative is a special case of the multivariate Riemann-Liouville derivative defined in \cite{malk2015}, which appears as a straightforward generalisation of the case for integer $r_i$:

\begin{equation*}
\frac{{{\partial ^{\sum {{r_i}} }}}}{{\partial {x_1}^{{r_1}} \ldots \partial {x_n}^{{r_n}}}}\prod\limits_i {x_i}^{p_i} = \prod\limits_i \frac{{\Gamma ({p_i} + 1)}}{{\Gamma ({p_i} - {r_i} + 1)}}{x_1}^{{p_i} - {r_i}}
\end{equation*}

The expression for the multivariate Gr\"unwald-Letnikov derivative (with $h$ and $N$ allowed to vary), as per Eq.~\eqref{mvglder-def}, can be factorised as follows:

\begin{equation*}
\frac{{{\partial ^{\sum {{r_i}} }}}}{{\partial {x_1}^{{r_1}} \ldots \partial {x_n}^{{r_n}}}}\prod\limits_i {x_i}^{p_i} = \prod\limits_i
\left[{h_i}^{-r_i}\sum\limits_{k_i = 0}^{N_i} (-1)^{k_i}\binom{r_i}{k_i}(x_i - k_i h_i) ^ {p_i}\right]
\end{equation*}

Where we suppress the limits $h_i \to 0$ and $N_i \to \infty$. The key difference between this and the univariate case is that we now have $n$ pairs $(N_i, h_i)$, which can in principle interrelate. The fact that both terms are fully factorised means that we can write down the characteristic equation corresponding to the multivariate case as precisely the \textit{product} of the expressions in the component-wise characteristic equations for the univariate case from \cite{pasu2018}, setting $q_i = N_i h_i/x_i$:

\begin{equation}
\prod\limits_i {{q_i}^{ - {r_i}}\sum\limits_{j = 0}^{{p_i}} {\binom{p_i}{j}} \frac{{{{( - {q_i})}^j}}}{{{r_i} - j}}}  = \prod\limits_i {\frac{\pi }{{\sin \pi {r_i}}}\binom{p_i}{r_i}} 
\label{drxp-char}
\end{equation}

Or alternatively:

\begin{equation}
\prod\limits_i \left[\frac{\pi r_i}{\sin\pi r_i}\binom{p_i}{r_i}{q_i}^{r_i} \frac{1}{{}_2F_1(- p_i, - r_i; 1 - r_i; q_i)}\right] = 1
\label{drxp-char2}
\end{equation}

One may check that $q_i = 1$ is always a solution to this equation.

\section{Conclusion}\label{sec:5}

\setcounter{section}{5}
\setcounter{equation}{0}\setcounter{theorem}{0}

Our formalism of infinitesimal translations has allowed us to concisely define the multivariate Gr\"unwald-Letnikov derivative in a well-motivated way. The result is a clear generalisation of the well-known expression for the univariate derivative.

In addition, we have found the general characteristic equation for the principal value of any derivative of a multivariate power function --  interestingly, the solution $q_i = 1$ always yields the principal value, implying that it works for any analytic function (as their Taylor expansion expresses them as an infinite sum of power functions). These results would be more significant under a proof of equivalence between the multivariate Gr\"unwald-Letnikov and Riemann-Liouville derivatives for functions where both derivatives are defined, analogous to the proof in \cite{orti2004} for the univariate case -- this problem will be the focus of future work.

It is interesting to consider the generalisation of the rank-$r$ gradient tensor (i.e. the tensor $\partial_{i_1} \ldots \partial_{i_r}f$). For the integer-order partial derivatives, the rank-$r$ gradient tensor of a function $f : \mathbb{R}^n \to \mathbb{R}$ allows for co-ordinates from ${\left[ {1,n} \right]^r}$ -- this follows from having $n$ choices for each successive derivative acting on the function. By contrast for the fractional-order derivative, there is a continuum of successive derivatives acting on the function, and the co-ordinate realization of the gradient tensor with ``real rank'' is the set of \textit{functions} $\{\tau : \left[0, r\right]_{\mathbb{R}} \to \left[1, n\right]_{\mathbb{Z}} \}$ which we'll call ``co-ordination functions'' and the corresponding generalisation of the symmetry of mixed partial derivatives is as follows: 

\begin{equation}
{D_\tau } = {D_\sigma } \Leftrightarrow \forall i \in \left[ {1,n} \right],\mu \left( {{\tau ^{ - 1}}(i)} \right) = \mu \left( {{\sigma ^{ - 1}}(i)} \right)
\label{fract-symm}
\end{equation}

Where $D_\tau$ is the derivative operator determined by the co-ordination function $\tau$, $\tau^{-1}(i)$ is the pre-image set of $i$ under $\tau$ and $\mu(S)$ is the measure of  $S \subseteq \mathbb{R}$.

This understanding of the real-rank tensors might, if explored further, shed some light on results relating to the the transformation laws for the Riemann-Liouville derivative derived in \cite{malk2015}.




 \bigskip \smallskip

 \it

 \noindent
$^1$ Department of Mathematics \\
Imperial College of London \\
180 Queen's Gate, South Kensington Campus \\
London -- SW7 2AZ, UNITED KINGDOM  \\[4pt]
  e-mail: ap6218@imperial.ac.uk
\hfill Received: April 4, 2019 \\[12pt]
\vspace*{-160pt}
\end{document}